# Asymptotics and Numerics of Zeros of Polynomials That Are Related to Daubechies Wavelets

Nico M. Temme

*CWI, P.O. Box 94079, 1090 GB Amsterdam, The Netherlands*

*e-mail:* `nicot@cwi.nl`

## Abstract

We give asymptotic approximations of the zeros of certain high degree polynomials. The zeros can be used to compute the filter coefficients in the dilation equations which define the compactly supported orthogonal Daubechies wavelets. Computational schemes are presented to obtain the numerical values of the zeros within high precision.

*1991 Mathematics Subject Classification:* 33B20, 41A60, 65D20.

## 1. Daubechies wavelets

The polynomial formed by the first $N$ terms of the binomial expansion

$$(1-y)^{-N} = \sum_{k=0}^{\infty} \binom{k+N-1}{k} y^k = \sum_{k=0}^{\infty} \frac{(N)_k}{k!} y^k$$

that is, the polynomial

$$P_N(y) = \sum_{k=0}^{N-1} \binom{k+N-1}{k} y^k = 1 + Ny + \frac{N(N+1)}{2} y^2 + \ldots + \binom{2N-2}{N-1} y^{N-1}, \ (1.1)$$

where $(a)_k = \Gamma(a+k)/\Gamma(a)$, plays an important role in the construction of the compactly supported Daubechies wavelets. There is a close connection between the zeros of $P_N(y)$ and the $2N$ filter coefficients $h(n)$ of the Daubechies wavelet $D_{2N}$. For a complete account of the theory we refer to Chapter 6 of DAUBECHIES (1992). In this section we give the details that are relevant to this paper. In later sections we describe the asymptotic methods for obtaining the zeros of $P_N(y)$ and discuss methods how to obtain the coefficients $h(n)$ numerically. This paper is an extension of work done by SHEN & STRANG (1995).

## 1.1. Some properties of Daubechies wavelets

We recall that the filter coefficients $h(n)$ define the dilation equation

$$\phi(x) = \sqrt{2} \sum_{n=0}^{2N-1} h(n)\phi(2x-n), \qquad (1.2)$$

the solution of which is called the scaling function. We take the following normalization of $\phi$ and, hence, of the coefficients $h(n)$:

$$\int_{-\infty}^{\infty} \phi(x)\,dx = 1, \qquad \sum_{n=0}^{2N-1} h(n) = \sqrt{2}. \qquad (1.3)$$

For the Daubechies wavelets the filter coefficients are real.

When the filter coefficients and the scaling function $\phi$ are available, the corresponding compactly supported orthogonal Daubechies wavelet is given by

$$\psi(x) = \sqrt{2} \sum_{n=2-2N}^{1} (-1)^n h(1-n)\phi(2x-n), \qquad (1.4)$$

which is also denoted by $D_{2N}$.

Denote the Fourier transform of $\phi$ by

$$\widehat{\phi}(\xi) = \frac{1}{\sqrt{2\pi}} \int_{-\infty}^{\infty} e^{-ix\xi}\phi(x)\,dx. \qquad (1.5)$$

Then (1.2) and (1.4) give the relations

$$\widehat{\phi}(\xi) = m_0(\xi/2)\,\widehat{\phi}(\xi/2), \qquad \widehat{\psi}(\xi) = \mu_0(\xi/2)\,\widehat{\phi}(\xi/2), \qquad (1.6)$$

where

$$\begin{aligned}
m_0(\xi) &= \frac{1}{\sqrt{2}} \sum_{n=0}^{2N-1} h(n)e^{-in\xi}, \\
\mu_0(\xi) &= \frac{1}{\sqrt{2}} \sum_{n=2-2N}^{1} (-1)^n h(1-n)e^{-in\xi} \\
&= -e^{-i\xi} m_0(-\xi-\pi).
\end{aligned} \qquad (1.7)$$

On the one hand, the orthonormality of the functions $\phi(\cdot + k), k \in \mathbb{Z}$ can be expressed in terms of the filter coefficients by the relation

$$\sum h(n-2j)h(n-2k) = \delta_{j,k},$$

where $j, k, n$ run over all relevant indices that define the finite set $h(n)$. On the other hand, the orthonormality is given by the relation

$$|m_0(\xi)|^2 + |m_0(\xi+\pi)|^2 = 1. \qquad (1.8)$$

The trigonometric polynomial $m_0(\xi)$ plays an important role in wavelet theory. It is called the filter function or transfer function. In classifying the compactly supported orthogonal wavelets, Daubechies assumes that the $m_0(\xi)$ has an $N$−fold zero at $\pm\pi$. That is, she writes

$$m_0(\xi) = \left(\frac{1 + e^{-i\xi}}{2}\right)^N \sum_{n=0}^{N-1} f(n)e^{-in\xi}. \tag{1.9}$$

A consequence is the $N$−fold zero of $m_0(\xi + \pi)$ and $m_0(-\xi - \pi)$ at $\xi = 0$. Using the relation for $\widehat{\psi}$ in (1.6) and $\mu_0$ in (1.7) we see that the $N$−fold zero of $m_0(\xi)$ at $-\pi$ is related to the vanishing of the $N$ moments

$$\int \psi(x)x^k \, dx, \quad k = 0, 1 \ldots, N-1.$$

This property gives a certain degree of approximating quality of the wavelet and leads to vanishing moments for the filter coefficients:

$$\sum_{n=0}^{2N-1} (-1)^n n^k h(n) = 0, \quad k = 0, 1, \ldots N-1. \tag{1.10}$$

## 1.2. The construction of Daubechies wavelets

Taking the symmetric product of the two filter functions, and writing $z = e^{i\xi}$, we have

$$m_0(\xi)m_0(-\xi) = |m_0(\xi)|^2 = \cos^{2N}\left(\tfrac{1}{2}\xi\right) P_N(y), \tag{1.11}$$

where $P_N(y)$ is the polynomial defined by

$$P_N(y) = Q_N(z)Q_N(1/z), \quad Q_N(z) = \sum_{n=0}^{N-1} f(n)z^{-n} \tag{1.12}$$

and $y$ is defined by $z + 1/z = 2 - 4y$, that is, $y = \sin^2 \tfrac{1}{2}\xi$. Indeed, the product $Q_N(z)Q_N(1/z)$ can be written as a polynomial of $y$ or $\cos\xi$, as can be verified easily.

When the coefficients $f(n)$ in (1.12) are known, the filter coefficients $h(n)$ follow from (1.9) and (1.7). The problem for constructing Daubechies wavelets is to find $f(n)$.

Substituting (1.11) into the orthonormality relation (1.8) we obtain for the polynomial $P_N$ the relation

$$(1 - y)^N P_N(y) + y^N P_N(1 - y) = 1. \tag{1.13}$$

A polynomial solution of this equation is the polynomial introduced in (1.1). To see this we introduce the incomplete beta function

$$I_x(a, b) = \frac{1}{B(a, b)} \int_0^x t^{a-1}(1 - t)^{b-1} \, dt, \quad B(a, b) = \frac{\Gamma(a)\Gamma(b)}{\Gamma(a + b)}, \tag{1.14}$$

where $\Re a > 0, \Re b > 0$ and $B(a, b)$ is the (complete) beta integral. We have $I_x(a, b) + I_{1-x}(b, a) = 1$ and we verify easily that the function

$$P_N(y) = (1 - y)^{-N} I_{1-y}(N, N) \tag{1.15}$$

satisfies the relation in (1.13). Also, this $P_N$ is the polynomial introduced in (1.1). Namely, using the substitution $t = (1 - y)u$ we obtain from (1.14):

$$
\begin{aligned}
(1 - y)^{-N} I_{1-y}(N, N) &= \frac{(1 - y)^{-N}}{B(N, N)} \int_0^{1-y} t^{N-1}(1 - t)^{N-1} \, dt \\
&= \frac{1}{B(N, N)} \int_0^1 u^{N-1}(1 - u + yu)^{N-1} \, du \\
&= \frac{1}{B(N, N)} \sum_{k=0}^{N-1} \binom{N-1}{k} y^k \int_0^1 u^{N-1+k}(1 - u)^{N-1-k} \, du \\
&= \sum_{k=0}^{N-1} \binom{N-1}{k} y^k \frac{(N+k-1)!(N-k-1)!}{(N-1)!(N-1)!} \\
&= \sum_{k=0}^{N-1} \binom{N+k-1}{k} y^k.
\end{aligned}
$$

## 1.3. Determining the filter coefficients

When we know the quantity $P_N$, we can try to solve for $Q_N$ in (1.12), and/or to find the $f(n)$. As an example consider $N = 2$. From (1.12) we obtain

$$
\begin{aligned}
1 + 2y &= [f(0) + f(1)/z][f(0) + f(1)z] \\
&= f(0)^2 + f(1)^2 + f(0)f(1)(z + 1/z) \\
&= ([f(0) + f(1)]^2 - 4yf(0)f(1),
\end{aligned}
$$

from which we can solve for $f(0), f(1)$. From (1.3), (1.7) and (1.9) we already know that $f(0) + f(1) = 1$. Furthermore we obtain $[f(0) - f(1)]^2 = 3$. This gives $f(0) = (1 + \sqrt{3})/2, f(1) = (1 - \sqrt{3})/2$. The filter coefficients follow from the relation (see ((1.7) and (1.9))

$$[(1 + 1/z)/2]^2 [f(0) + f(1)/z] = \frac{1}{\sqrt{2}}[h(0) + h(1)/z + h(2)/z^2 + h(3)/z^3].$$

This gives the coefficients of $D_4$:

$$h(0) = \frac{1 + \sqrt{3}}{4\sqrt{2}}, \quad h(1) = \frac{3 + \sqrt{3}}{4\sqrt{2}}, \quad h(2) = \frac{3 - \sqrt{3}}{4\sqrt{2}}, \quad h(3) = \frac{1 - \sqrt{3}}{4\sqrt{2}}. \tag{1.16}$$

For $N = 3$ we can also obtain exact solutions of the nonlinear equations, but for larger values of $N$ this is not possible. In addition, the complexity of the computational scheme increases.

**Remark 1.1.** The above equations for $f(0), f(1)$ are symmetric with respect to these quantities. Interchanging the values of $f(0), f(1)$ gives a different set $h(n)$ (with $h(j) \to h(3-j), j = 0, 1, 2, 3$). The present choice of $f(0), f(1)$ gives $Q_2(z) = f(0) + f(1)/z$, with a zero at $z_1 = 2 - \sqrt{3}$, inside the unit circle.

A different method is based on the zeros of $P_N(y)$. When we have computed the $N-1$ zeros $y_n$ of $P_N(y)$, we can compute the corresponding zeros $z_n$ of $Q_N(z)$, by using the relation $z_n + 1/z_n = 2 - 4y_n$. We have (recall that $Q_N(1) = 1$)

$$P_N(y) = Q_N(z)Q_N(1/z) = \prod_{n=1}^{N-1} \frac{1 - z_n/z}{1 - z_n} \frac{1 - zz_n}{1 - z_n}$$
$$= \prod_{n=1}^{N-1} \frac{z_n(z_n + 1/z_n - z - 1/z)}{(1 - z_n)^2}$$
$$= (4z_n)^{N-1} \prod_{n=1}^{N-1} \frac{y - y_n}{(1 - z_n)^2}.$$

For each value $y_n$ there are two $z-$zeros: $z_n$ and $1/z_n$. We use the zeros $z_n$ inside the unit circle. By expanding the product representation of $Q_N$ we obtain the coefficients $f(n)$ of (1.12), which are needed in (1.9). By expanding and using (1.7) and (1.9) we obtain the coefficients $h(n)$ from

$$\frac{1}{\sqrt{2}} \sum_{n=0}^{2N-1} h(n)z^{-n} = \left(\frac{1 + 1/z}{2}\right)^N \sum_{n=0}^{N-1} f(n)z^{-n}. \tag{1.17}$$

We verify this second method by taking $N = 2$. We have $P_2(y) = 1 + 2y$, with zero $y_1 = -\frac{1}{2}$. The corresponding $z_1$ inside the unit circle is $z_1 = 2 - \sqrt{3}$. Equation (1.17) reads

$$\frac{1}{\sqrt{2}} \sum_{n=0}^{3} h(n)z^{-n} = \left(\frac{1 + 1/z}{2}\right)^2 \frac{1 - z_1/z}{1 - z_1}.$$

This gives the coefficients $h(n)$ given in (1.16).

In Figure 1.1 we show the 99 zeros of $P_{100}(y)$ and the corresponding $z-$zeros.

In Shen & Strang (1995) the location of the zeros of $P_N(y)$ for large values of $N$ is discussed. The zeros approach a limiting curve $|4y(1-y)| = 1$ in the complex plane. As remarked by Shen and Strang, the wide dynamic range in the coefficients of $P_N(y)$ makes the zeros difficult to compute for large $N$. They give pictures of the zeros located near the limiting curve up to $N = 70$, and they observe false results when using standard Matlab methods for $N = 100$. In the present paper we give an asymptotic expansion of the zeros, and we show by numerical verification that these approximations are excellent starting values for obtaining iterated high precision values of the zeros. The methods of computing the high degree polynomial $P_N(y)$ are also based on a uniform asymptotic representation of the incomplete beta function, of which details are given in a final section.

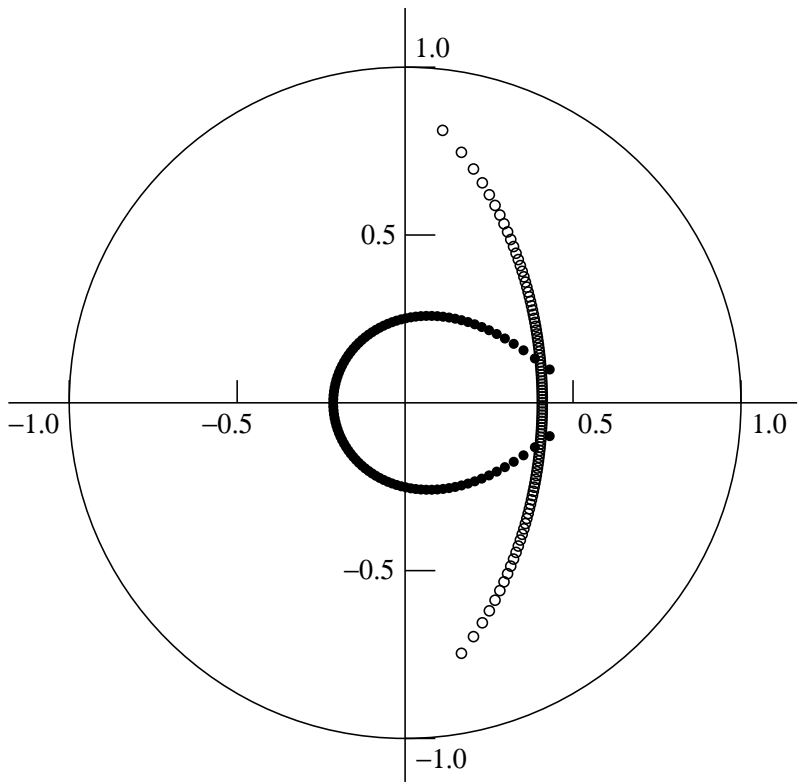

**Figure 1.1.** Zeros $y_n$ (black dots) of $P_N(y)$ for $N = 100$ and corresponding zeros $z_n$ (open dots) defined by $z_n + 1/z_n = 2 - 4y_n$.

## 2. Asymptotic inversion of the incomplete beta function

We recall the relation given in (1.15)

$$P_N(y) = (1 - y)^{-p} I_{1-y}(p, p), \tag{2.1}$$

where the incomplete beta function is defined in (1.14). In TEMME (1992) the asymptotic inversion of the incomplete beta function $I_x(a, b)$ for large values of $a$ and $b$ is discussed. That is, an asymptotic approximation is given of the solution $x$ of the equation $I_x(a, b) = q$, with $0 < q < 1$. We can use the same inversion method for finding complex $x$−zeros (different from zero) by taking $q = 0$. In TEMME (1992) three cases have been considered:

- $a$ and $b$ large, with $a - b$ bounded;
- $a$ and $b$ large, with $a/b$ and $b/a$ bounded away from zero;
- at least one of the parameters $a, b$ is large.

The first case applies to the present situation, because we are interested in the zeros of $I_{1-y}(N, N)$.

We repeat the main steps of the inversion procedure of $I_x(a, a + \beta)$ given in §2 of TEMME (1992), by taking $a = N, \beta = 0$. First we write $I_{1-y}(N, N)$ in terms of an error function.

## 2.1. Transformation into Gaussian form

We have

$$I_{1-y}(N, N) = \frac{4^{-N}}{B(N, N)} \int_0^{1-y} [4t(1-t)]^N \frac{dt}{t(1-t)}.$$

We transform this to a standard form with a Gaussian character by writing

$$\begin{aligned}
-\tfrac{1}{2}\zeta^2 &= \ln[4t(1-t)], \quad 0 < t < 1, \quad \text{sign}(\zeta) = \text{sign}(t - \tfrac{1}{2}), \\
-\tfrac{1}{2}\eta^2 &= \ln[4y(1-y)], \quad 0 < y < 1, \quad \text{sign}(\eta) = \text{sign}(\tfrac{1}{2} - y).
\end{aligned} \tag{2.2}$$

Therefore

$$I_{1-y}(N, N) = \frac{4^{-N}}{B(N, N)} \int_{-\infty}^{\eta} e^{-\frac{1}{2}N\zeta^2} \frac{1}{t(1-t)} \frac{dt}{d\zeta} d\zeta.$$

We can invert the relations in (2.2):

$$\begin{aligned}
t &= \tfrac{1}{2}\left[1 \pm \sqrt{1 - \exp(-\tfrac{1}{2}\zeta^2)}\right] = \tfrac{1}{2}\left[1 + \zeta \sqrt{[1 - \exp(-\tfrac{1}{2}\zeta^2)]/\zeta^2}\right], \\
y &= \tfrac{1}{2}\left[1 \pm \sqrt{1 - \exp(-\tfrac{1}{2}\eta^2)}\right] = \tfrac{1}{2}\left[1 - \eta \sqrt{[1 - \exp(-\tfrac{1}{2}\eta^2)]/\eta^2}\right],
\end{aligned} \tag{2.3}$$

where the second square roots are non-negative for real values of their arguments. It easily follows that

$$\frac{1}{t(1-t)} \frac{dt}{d\zeta} = \frac{-\zeta}{1 - 2t},$$

and that the following standard form (in the sense of TEMME (1982)) can be obtained

$$I_{1-y}(N, N) = \sqrt{\frac{N}{2\pi}} \, \Phi(N) \int_{-\infty}^{\eta} e^{-\frac{1}{2}N\zeta^2} \phi(\zeta) \, d\zeta, \tag{2.4}$$

where

$$\Phi(N) = \frac{1}{\sqrt{N}} \frac{\Gamma(N + \tfrac{1}{2})}{\Gamma(N)}, \quad \phi(\zeta) = \sqrt{\frac{\tfrac{1}{2}\zeta^2}{1 - \exp(-\tfrac{1}{2}\zeta^2)}}. \tag{2.5}$$

We have

$$\Phi(N) \sim 1 - \tfrac{1}{8}N^{-1} + \tfrac{1}{128}N^{-2} + \ldots \qquad (N \to \infty). \tag{2.6}$$

The function $\phi(\zeta)$ is analytic in a strip containing $\mathbb{R}$; the singularities nearest to the origin occur at $\pm 2\sqrt{\pi} \exp(\pm i\pi/4)$. The first part of the Taylor expansion of this even function reads:

$$\phi(\zeta) = 1 + \frac{1}{8}\zeta^2 + \frac{1}{384}\zeta^4 + \ldots$$

We write (2.4) in the form (see TEMME (1982))

$$I_{1-y}(N, N) = \tfrac{1}{2}\text{erfc}\left(-\eta\sqrt{N/2}\right) + R_N(\eta), \tag{2.7}$$

where $\eta$ is defined in (2.2), and erfc is the error function defined by

$$\text{erfc } z = \frac{2}{\sqrt{\pi}} \int_z^{\infty} e^{-t^2} dt. \tag{2.8}$$

In fact, we replace in (2.4) $\Phi(N)\phi(\zeta)$ by 1, and the error is contained in $R_N(\eta)$. We try to find zeros of $I_{1-y}(N, N)$ by using representation (2.7), assuming that $N$ is large. First we find zeros in terms of $\eta$; afterwards, we determine $y$ from the second line in (2.3). When $N$ is large, we consider the zeros of the error function as a first approximation to the zeros of the incomplete beta function.

## 2.2. Inversion of the incomplete beta function

Let $\eta_0, \eta$ solve the equations

$$\frac{1}{2}\operatorname{erfc}\left(-\eta_0\sqrt{\tfrac{1}{2}N}\right) = 0,$$
$$\frac{1}{2}\operatorname{erfc}\left(-\eta\sqrt{\tfrac{1}{2}N}\right) + R_N(\eta) = 0. \tag{2.9}$$

We assume that $\eta_0$ is known and that, for large $N$, this value is a first approximation for the value $\eta$ that satisfies the second equation in (2.9). The relation between $\eta$ and $\eta_0$ is written as

$$\eta = \eta_0 + \varepsilon, \tag{2.10}$$

and we try to determine $\varepsilon$. We can expand the quantity $\varepsilon$ in the form

$$\varepsilon \sim \frac{\varepsilon_1}{N} + \frac{\varepsilon_2}{N^2} + \frac{\varepsilon_3}{N^3} + \ldots, \tag{2.11}$$

as $N \to \infty$. The coefficients $\varepsilon_i$ can be found by using a perturbation method. For details we refer to TEMME (1992).

The first coefficient $\varepsilon_1$ is given by

$$\varepsilon_1 = \frac{1}{\eta_0}\ln\phi(\eta_0), \tag{2.12}$$

where $\phi$ is defined in (2.5). The quantity $\varepsilon_1$ is an odd function of $\eta_0$ and analytic in the disc $|\eta_0| < 2\sqrt{\pi}$. Further terms $\varepsilon_i$ are given by

$$\varepsilon_2 = \frac{1}{8\eta\phi}(8\phi\varepsilon_1' + 8\phi'\varepsilon_1 - \phi - 4\phi\varepsilon_1^2); \tag{2.13}$$

$$\varepsilon_3 = \frac{1}{128\eta\phi}(128\phi\varepsilon_2' + 128\phi'\varepsilon_1\varepsilon_1' - 16\phi\varepsilon_1' + 128\phi'\varepsilon_2 + 64\phi''\varepsilon_1^2 + \\ -\,128\phi'\varepsilon_1 + \phi - 128\phi\varepsilon_1\varepsilon_2 - 64\phi\varepsilon_2^2\eta^2 - 64\phi\varepsilon_2\eta\varepsilon_1^2 - 16\phi\varepsilon_1^4). \tag{2.14}$$

The derivatives $\phi', \varepsilon'$, etc., are with respect to $\eta$, and all quantities are evaluated at $\eta_0$. More information on the coefficients $\varepsilon_i$ will be given in §4.2.

## 3. More details on the zeros and numerical examples

As explained in the previous section, the zeros of $P_N(y)$ are approximated in terms of the zeros of the error function. From FETTIS, CASLIN & CRAMER (1973) we know that two infinite strings of zeros of erfc $w$ occur in the neighbourhood of the diagonals $v = \pm u$ in the left half plane $u < 0$, $w = u + iv$. The first few zeros are given in Table 3.1. Numerical values of the first 100 zeros of erfc $w$ and asymptotic approximations of the zeros are also given by Fettis et al. A first order approximation reads

$$w_k^\pm \sim \sqrt{2\pi(k-1/8)}\ e^{\pm 3\pi i/4}, \quad k \to \infty. \tag{3.1}$$

| $k$ | $u_k$ | $\pm$ | $i\,v_k$ |
|---|---|---|---|
| 1 | $-1.35481...$ | $\pm i$ | $1.99146...$ |
| 2 | $-2.17704...$ | $\pm i$ | $2.69114...$ |
| 3 | $-2.78438...$ | $\pm i$ | $3.23533...$ |
| 4 | $-3.28741...$ | $\pm i$ | $3.69730...$ |
| 5 | $-3.72594...$ | $\pm i$ | $4.10610...$ |

**Table 3.1.** First five pairs $w_k^\pm = u_k \pm i\,v_k$ of zeros of erfc $w$

When $\eta$ is a zero of the right-hand side of (2.7), we use (2.10)–(2.11) for obtaining an asymptotic expansion of $\eta$ for large positive values of $N$; $-\eta_0\sqrt{N/2}$ is a zero of the error function. The mapping of the $\eta$−plane to the $y$−plane is given by (see (2.2) – (2.3))

$$-\tfrac{1}{2}\eta^2 = \ln\left[4y(1-y)\right]. \tag{3.2}$$

Because the zeros in terms of $\eta$ occur in the neighbourhood of the diagonals $\Re\eta = \pm\Im\eta$ in the right half of the $\eta$−plane, we see that the $y$−zeros of $P_N(y)$ occur near the curve defined by $|4y(1-y)| = 1$, with $\Re y < \tfrac{1}{2}$. The full curve is a lemniscate; the extreme points cut the real axis at $y = \tfrac{1}{2} \pm \tfrac{1}{\sqrt{2}}$.

The mapping in (3.2) is singular at the points

$$\pm 2\sqrt{k\pi}\ e^{\pm\frac{1}{4}\pi i}, \quad k = 1, 2, 3, \dots$$

on the diagonals in the $\eta$−plane. All these branch points are mapped to $y = \tfrac{1}{2}$. Of special interest are the points on the diagonals in the right half plane, for instance the points $2\sqrt{\pi}\exp(\pm\pi i/4)$. The parts of the diagonals given by $\eta = \rho e^{\pm\pi i/4}, 0 \le \rho \le 2\sqrt{\pi}$ are mapped to the left leave of the lemniscate, the parts satisfying $2\sqrt{\pi} \le \rho \le 2\sqrt{2\pi}$, are mapped to the right leave, and so on.

If $N$ is even $P_N(y)$ has one real zero (recall that this polynomial is of degree $N-1$), and this real zero will be approximated by using $w_{N/2}^\pm$ (see (3.1)), and for the complex zeros we can use $w_k^\pm, k = 1, 2, \dots N/2 - 1$. When $N$ is odd we use $w_k^\pm, k = 1, 2, \dots (N-1)/2$.

We have used the asymptotic expansion (2.11) for $N = 100$ with five terms, and obtained approximations of the zeros of $P_N(y)$ which were accurate with a precision of at least 12 digits. The zeros near $y = \tfrac{1}{2}$ are more accurate than those near the extreme point $y = \tfrac{1}{2} - \tfrac{1}{\sqrt{2}}$. The latter correspond to $\eta$−values that are closer to the singular points $2\sqrt{\pi}\exp(\pm\pi i/4)$ of the mapping given in (2.2). The asymptotic expansion (2.11) breaks down at these points.

When $N$ increases the used $\eta$−zeros are bounded away from these singular points. To see this, observe that we need $w_k^\pm$, with $1 \le k \le N/2$. From (3.1) we see that $w_{N/2}^\pm \sim$ $\sqrt{N\pi}\,e^{\pm 3\pi i/4}$. Because $w_k$, a zero of the error function, corresponds to $-\eta_0\sqrt{N/2}$, we see that the maximal $\eta_0$ satisfies $\eta_0 \sim \sqrt{2\pi}\,\exp(\pm\pi i/4)$.

**Remark 3.1.** We will not discuss the role of other zeros in the $\eta-$plane in the neighborhood of the diagonals. A similar phenomenon is explained in TEMME (1995). The point is that, although we started with integer values of $N$, the incomplete beta function in (2.1) and (2.7) can be interpreted with general complex values of $N$. When $N$ is not an integer there are more $y-$zeros, but those occur on other Riemann sheets of the multi-valued incomplete beta function.

### 3.1. Numerical aspects

When $N$ is large the representation of $P_N(y)$ in (1.1) is not suitable for verifying the accuracy of the zeros that are obtained by asymptotic methods with, say, 12 accurate digits. When $N = 100$ the coefficient of the power $y^{99}$ is an integer with about 60 digits, and $\binom{198}{99}y^{99}$ is about $10^{30}$ when $y = \frac{1}{2}$. It follows that the zeros have to be very accurate to verify them by using (1.1), and we have to perform the computations with many extra digits.

The incomplete beta function can computed for large values of $N$ by using a uniform asymptotic expansion, which a gives a more stable representation. More details are given in the next section. We have verified the accuracy of the zeros by using this expansion, and obtained the accuracy mentioned above. However, we can use different numerical verifications, for instance the ones based on sums and products of zeros of polynomials. For the $n$ zeros $y_k$ of the polynomial $P(y) = a_n y^n + \ldots + a_0$, we have the relations

$$a_n \sum_{k=1}^{n} y_k = -a_{n-1}, \quad a_n \prod_{k=1}^{n} y_k = (-1)^n a_0.$$

In the present case we know that the zeros $y_k$ of $P_N(y)$ have to satisfy the rules

$$\sum_{k=1}^{n} y_k = -\frac{1}{2}, \quad \binom{2N-2}{N-1}\prod_{k=1}^{N-1} y_k = (-1)^{N-1}.$$

We verified both rules for $N = 100$ and we obtained by using 5 terms in (2.11)

$$\sum_{k=1}^{99} y_k + \frac{1}{2} = -2.56\cdots \times 10^{-13}, \quad \binom{198}{99}\prod_{k=1}^{99} y_k + 1 = -1.56\cdots \times 10^{-12},$$

which confirms our earlier claim about the accuracy of the zeros.

We have obtained more accurate values of the zeros of $P_N(y)$, $N = 100$ by using the multi-length facilities of Maple and using the earlier obtained zeros as starting values in a Newton-Raphson process based on representations of $P_N(y)$ given in (2.1) and (2.7). In this way we obtained the zeros $y_n$ with 45 relevant digits.

In SHEN & STRANG (1995) there is no discussion on numerical methods for obtaining the filter coefficients $h(n)$, once the zeros $y_n$ or $z_n$ are available. As we have explained in §1.3 we need the $f(n)$ in (1.17), but the computation of $f(n)$ by using the zeros of $Q_N(z)$ may be unstable when $N$ is large. Another source of errors comes from expanding the product of the binomial term and $Q_N(z)$ in (1.17). Unfortunately, we don't have available asymptotic methods for obtaining the filter coefficients (we have powerful methods for $P_N$, but not for $Q_N$). Further numerical research is needed to compute $h(n)$.

In fact, highly accurate values of $y_n$ are needed in order to obtain accurate values of $h(n)$ by using straightforward methods. Expanding the right-hand side of (1.17) in powers of $1/z$, and using 45-digits accurate values of $z_n$, we have obtained a set $h(n)$ that satisfies

$$\sum_{n=0}^{2N-1} h^2(n) - 1 = 5.42 \times 10^{-19}, \quad N = 100,\tag{3.3}$$

which shows that the dominant $h(n)$ are accurate with about 20 digits accuracy. A test based on the vanishing moments given in (1.10) indicated that only for rather low $k-$values a precision as in (3.3) could be obtained.

In Figure 3.1 we give an idea of the values of the filter coefficients; we give the values of $-^{10}\log|h(n)|, n = 0, 1, \ldots, 199$. The largest value of the coefficients is $h(20) = 0.39910\cdots$.

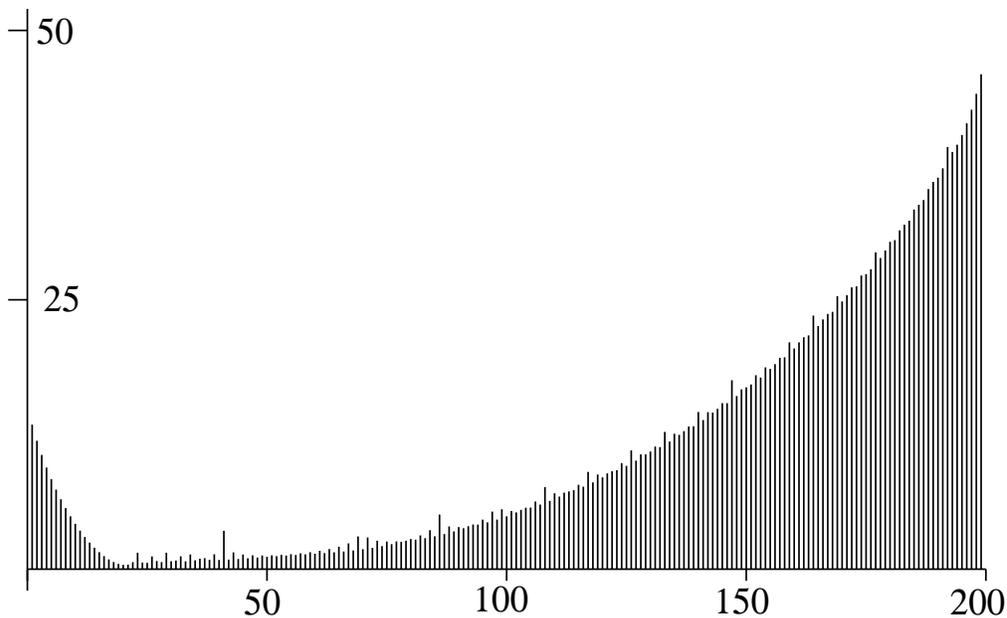

**Figure 3.1.** Filter coefficients $h(n)$ of the Daubechies wavelet $D_{200}$; shown are the values of $-^{10}\log|h(n)|, n = 0, 1, \ldots, 199$.

# 4. Asymptotic expansions

## 4.1. The incomplete beta function

We summarize from TEMME (1982), (1996) the details of an expansion for the incomplete beta function that is based on formula (2.7). We can expand

$$R_N(\eta) \sim \frac{e^{-\frac{1}{2}N\eta^2}}{\sqrt{2\pi N}} \sum_{k=0}^{\infty} \frac{B_k(\eta)}{N^k}, \quad N \to \infty. \tag{4.1}$$

This expansion holds uniformly in a strip $|\Im\eta| \le \sqrt{2\pi} - \delta$, where $\delta$ is a small positive number. The width of the strip is determined by the singularities of the mapping given in (2.2). As mentioned earlier, the singularities are defined by the $\eta-$values for which $\exp(-\frac{1}{2}\eta^2) = 1$. The $\eta-$zeros needed for $P_N(y)$ all lie within this strip.

By differentiating (2.4) and (2.7) with respect to $\eta$, substituting

$$\Phi(N) \sim \sum_{k=0}^{\infty} \frac{c_k}{N^k}$$

(a general form of (2.6)) and (4.1), we obtain a set of recursion relations for the coefficients $B_k$:

$$B_0 = \frac{1 - \phi(\eta)}{\eta}, \quad \eta B_{k+1} = \frac{d\, a\, set\, of\, B_k}{d\eta} - c_{k+1}\,\phi(\eta), \quad k = 0, 1, 2, \dots$$

with $\phi$ given in (2.5). The next few $B_k$ are

$$B_1 = -\frac{1}{8\eta^3}(3\phi\eta^2 - 8\phi^3 + 8),$$

$$B_2 = -\frac{1}{128\eta^5}(25\phi\eta^4 - 240\eta^2\phi^3 + 384\phi^5 - 384).$$

In these formulas we have removed the derivatives of $\phi$ by using

$$\frac{d\phi(\eta)}{d\eta} = \frac{\phi(\eta)}{\eta} \left[ 1 + \frac{1}{2}\eta^2 - \phi^2(\eta) \right], \tag{4.2}$$

which easily follows from (2.5).

For verifying the numerical computations of the zeros we have used expansion (4.1) with terms up to and including $B_{10}/N^{10}$.

## 4.2. The coefficients $\varepsilon_k$ of (2.11)

In (2.14) and (2.15) we have given the coefficients $\varepsilon_2$ and $\varepsilon_3$. It is convenient to remove derivatives by using (4.2). In this way we obtain

$$\varepsilon_2 = \frac{1}{8}(8 + 3\eta^2 - 8\phi^2 + 4\eta^3\varepsilon_1 - 8\phi^2\varepsilon_1\eta - 4\varepsilon_1^2\eta^2)/\eta^3,$$

$$\varepsilon_3 = \frac{1}{16}(-40 + 5\eta^4 + 56\phi^4 - 16\phi^2 + 8\eta^2 - 38\eta^2\phi^2 - 16\eta\varepsilon_1 + 16\phi^4\varepsilon_1^2\eta^2 - 32\phi^2\varepsilon_1\eta^3$$
$$+ 48\phi^4\varepsilon_1\eta - 6\eta^3\varepsilon_1 + 16\phi^2\varepsilon_1\eta + 4\eta^5\varepsilon_1 - 8\eta^4\varepsilon_1^2\phi^2 + 16\phi^2\varepsilon_1^2\eta^2 + 8\varepsilon_1^3\eta^3 - 8\varepsilon_1^2\eta^4)/\eta^5.$$

To avoid cancellations for small values of $|\eta|$ we can use the expansions

$$\varepsilon_1 = \frac{1}{8}\eta - \frac{1}{192}\eta^3 + \frac{1}{92160}\eta^7 - \frac{1}{23224320}\eta^{11} + \frac{1}{4954521600}\eta^{15} - \frac{1}{980995276800}\eta^{17} + \dots \ .$$

$$\varepsilon_2 = \frac{1}{128}\eta - \frac{5}{1536}\eta^3 + \frac{7}{40960}\eta^5 + \frac{1}{81920}\eta^7 - \frac{407}{371589120}\eta^9 + \dots,$$

$$\varepsilon_3 = \frac{-5}{1024}\eta - \frac{11}{24576}\eta^3 + \frac{63}{327680}\eta^5 - \frac{823}{165150720}\eta^7 + \dots,$$

$$\varepsilon_4 = \frac{-21}{32768}\eta + \frac{37}{65536}\eta^3 + \frac{179}{5242880}\eta^5 + \dots,$$

$$\varepsilon_5 = \frac{399}{262144}\eta + \frac{219}{2097152}\eta^3 + \dots.$$

All these expansions have radius of convergence $2\sqrt{\pi}$.